\documentclass[submission%
]{dmtcs}

\usepackage[latin1]{inputenc}
\usepackage{subfigure}

\usepackage[round]{natbib}

\usepackage{amsmath}
\usepackage{amssymb}
\usepackage{amsfonts}
\usepackage{tikz}
\usetikzlibrary{arrows,decorations.pathreplacing,patterns}
\usepackage{mathrsfs}

\def\Q{\check{Q}}
\def\wa{\widetilde{W}}

\def\Z{\mathbb{Z}}
\def\R{\mathbb{R}}
\def\P{\mathsf{Park}}
\def\D{\mathcal{D}}

\def\ac{A_{\circ}}

\def\wshi{W_{\mathsf{Shi}}}

\def\diag{\mathsf{Diag}}
\def\B{\mathcal B}

\newcommand{\sk}{

\smallskip}
\newcommand{\hypoct}{\mathfrak H}
\renewcommand{\S}{\mathfrak S}
\newcommand{\op}[1]{\operatorname{#1}}
\newcommand{\id}{\op{id}}
\newtheorem{mythrm}{Theorem}[section]
\newtheorem{mylem}[mythrm]{Lemma}
\newtheorem{myprop}[mythrm]{Proposition}
\newcommand{\reft}[1]{Theorem~\ref{thrm:#1}}
\newcommand{\refl}[1]{Lemma~\ref{thrm:#1}}

\newcommand{\reff}[1]{Figure~\ref{Figure:#1}}

\author{Robin Sulzgruber\thanks{Email: \email{robin.sulzgruber@univie.ac.at}.} \and Marko Thiel\thanks{Email: \email{marko.thiel@univie.ac.at}.\\Research supported by the Austrian Science Fund (FWF), grant S50-N15 in the framework of the Special Research Program ``Algorithmic and Enumerative Combinatorics'' (SFB F50).}}
\title{Type $C$ parking functions and a zeta map}
\address{Fakult\"at f\"ur Mathematik, Universit\"at Wien, Oskar-Morgenstern-Platz 1, 1090 Vienna, Austria}
\keywords{parking functions, Shi arrangement, zeta map, dinv statistic}
\received{2014-11-14}
\revised{tba}
\accepted{tba}

\begin{document}
\maketitle

\begin{abstract}
\paragraph{Abstract}
We introduce type $C$ parking functions, encoded as vertically labelled lattice paths and endowed with a statistic $\op{dinv'}$.
We define a bijection from type $C$ parking functions to regions of the Shi arrangement of type $C$, encoded as diagonally labelled ballot paths and endowed with a natural statistic $\op{area'}$.
This bijection is a natural analogue of the zeta map of Haglund and Loehr and maps $\op{dinv'}$ to $\op{area'}$.
We give three different descriptions of it. 

\paragraph{R\'esum\'e}
Nous introduisons les fonctions de stationnement de type $C$, encod\'ees par des chemins \'etiquet\'es verticalement et munies d'une statistique $\op{dinv'}$. Nous d\'efinissons une bijection entre les fonctions de stationnement de type $C$ et les r\'egions de l'arrangement de Shi de type $C$, encod\'ees par des chemins \'etiquet\'es diagonalement et munies d'une statistique naturelle $\op{area'}$. Cette bijection est un analogue naturel \`a la fonction zeta de Haglund et Loehr, et envoie $\op{dinv'}$ sur $\op{area'}$. Nous donnons trois diff\'erentes descriptions de celle-ci.
%
%
\end{abstract}

\section{Introduction and Motivation}\label{sec:intro}
One of the most well-studied objects in algebraic combinatorics is the space of diagonal harmonics of the symmetric group $\S_n$.
Its Hilbert series has two (conjectural) combinatorial interpretations: 
\begin{align*}
\mathcal{DH}(n;q,t)=\sum_{P\in\P_n}q^{\op{dinv'}(P)}t^{\op{area}(P)}=\sum_{R\in\diag_n}q^{\op{area'}(R)}t^{\op{bounce}(R)},
\end{align*}
where $\P_n$ is the set of \emph{parking functions} of length $n$, viewed as vertically labelled Dyck paths, and $\diag_n$ is the set of \emph{diagonally labelled} Dyck paths with $2n$ steps.
There is a bijection $\zeta$ due to 
\cite{HagLoe2005} that maps $\P_n$ to $\diag_n$ and sends the bistatistic $(\op{dinv'},\op{area})$ to $(\op{area'},\op{bounce})$, demonstrating the second equality.

\sk
The combinatorial objects $\P_n$ and $\diag_n$ may be viewed as the type $A_{n-1}$ cases of more general objects associated to any crystallographic root system $\Phi$.
These are, respectively, the \emph{finite torus} \mbox{$\Q/(h+1)\Q$} and the set of regions of the \emph{Shi arrangement} of $\Phi$. Here $\Q$ is the coroot lattice and $h$ is the Coxeter number of $\Phi$.
Both of these objects have the same cardinality $(h+1)^r$, where $r$ is the rank of $\Phi$, so there should be a uniform ``zeta map'' giving a bijection between them.
This map does in fact exist, and will be described in future work.

\sk
In the present extended abstract we focus on the root system of type $C_n$. In Section \ref{sec:defs} we present the necessary background on Weyl groups and the Shi arrangement. In Section \ref{sec:pf} we introduce combinatorial models for the finite torus of type $C$ in terms of \emph{vertically labelled lattice paths} and for the set of regions of the Shi arrangement of type $C$ in terms of \emph{diagonally labelled ballot paths}. We also introduce a statistic $\op{dinv'}$ on vertically labelled lattice paths and a statistic $\op{area'}$ on diagonally labelled ballot paths. These statistics are natural analogues of the corresponding statistics in type $A$.

In Section \ref{sec:zeta} we describe a map between these two combinatorial models that we call the \emph{type $C$ zeta map}. We give three descriptions of this map, all similar in style to different descriptions of the classical zeta map. The first description in terms of area vectors follows \cite{HagLoe2005}. The second description in terms of ascents and valleys resembles that of \cite[Section 5.2]{ALW2014}. The third description as a sweep map is in the spirit of \cite{ALW2014SweepMaps}. Our main result (\reft{HLzetaC}) is that the zeta map of type $C$ is a bijection that sends the $\op{dinv'}$ statistic to the $\op{area'}$ statistic.

\sk
We prioritise examples and prefer to include an adequate presentation of the known combinatorial objects of type $A$ rather than presenting proofs. A full version of this extended abstract containing all proofs is in preparation.

\section{Definitions and Preliminaries}\label{sec:defs}
\subsection{Weyl groups}

Let $\Phi$ be an irreducible crystallographic root system of rank $r$, with simple system $\Delta=\{\alpha_1,\alpha_2,\ldots,\alpha_r\}$, positive system $\Phi^+$ and ambient space $V$.
For background on root systems and reflection groups see \cite{Humph1990}.
For $\alpha\in\Phi$, let $s_{\alpha}$ be the reflection in the hyperplane 
\begin{align*}
H_{\alpha}=\{x\in V:\langle x,\alpha\rangle=0\}.
\end{align*}
Then the Weyl group $W$ of $\Phi$ is the group of automorphisms of $V$ generated by all the $s_{\alpha}$ with $\alpha\in\Phi$.
Define the \emph{Coxeter arrangement} of $\Phi$ as the central hyperplane arrangement in $V$ given by all the hyperplanes $H_{\alpha}$ for $\alpha\in\Phi$.
The connected components of the complement of the union of these hyperplanes are called \emph{chambers}.
The Weyl group $W$ acts simply transitively on the chambers, so if we define the \emph{dominant chamber} as
\begin{align*}
C=\{x\in V:\langle x,\alpha\rangle>0\text{ for all }\alpha\in\Delta\},
\end{align*}
we may write every chamber as $wC$ for a unique $w\in W$.

For $\alpha\in\Phi$ and $d\in\Z$, let $s_{\alpha}^d$ be the reflection in the affine hyperplane 
\begin{align*}
H_{\alpha}^d=\{x\in V:\langle x,\alpha\rangle=d\}.
\end{align*}
Then the affine Weyl group $\wa$ of $\Phi$ is the group of affine transformations of $V$ generated by all the $s_{\alpha}^d$ for $\alpha\in\Phi$ and $d\in\Z$.
Define the \emph{affine Coxeter arrangement} as the affine hyperplane arrangement in $V$ given by all the $H_{\alpha}^d$ for $\alpha\in\Phi$ and $d\in\Z$.
The connected components of the complement of the union of these hyperplanes are called \emph{alcoves}.
The affine Weyl group $\wa$ acts simply transitively on the alcoves, so if we write $\tilde{\alpha}$ for the highest root of $\Phi$ and define the \emph{fundamental alcove} as
\begin{align*}
\ac=\{x\in V:\langle x,\alpha\rangle>0\text{ for all }\alpha\in\Delta\text{ and }\langle x,\tilde{\alpha}\rangle<1\},
\end{align*}
we may write every alcove as $w_a\ac$ for a unique $w_a\in\wa$.
The affine Weyl group $\wa$ acts on the coroot lattice $\Q$, and if we identify $\Q$ with its translation group we may write $\wa=W\ltimes\Q$ as a semidirect product.

If $\alpha\in\Phi^+$ and $w_a\in\wa$, there is a unique integer $k$ such that $k<\langle x,\alpha\rangle<k+1$ for all $x\in w_aA_{\circ}$. We denote this integer by $k(w_a,\alpha)$.

\subsection{The Shi arrangement}\label{Subsection:shi}
Define the \emph{Shi arrangement} as the hyperplane arrangement given by the hyperplanes $H_{\alpha}^d$ for $\alpha\in\Phi^+$ and $d=0,1$. Then the complement of the union of these hyperplanes falls apart into connected components, which are called the \emph{regions} of the arrangement. 
The hyperplanes that support facets of a region $R$ are called the \emph{walls} of $R$. Those walls of $R$ that do not contain the origin and separate $R$ from the origin are called the \emph{floors} of $R$.
Define the walls and floors of an alcove similarly. Notice that every wall of a region is a hyperplane of the Shi arrangement, but the walls of an alcove need not be. We call a region or alcove \emph{dominant} if it is contained in the dominant chamber.
\begin{mythrm}\label{thrm:min} {\textnormal{\cite[Prop~7.1]{Shi1987}}}
Every region $R$ of the Shi arrangement has a unique \emph{minimal alcove} $w_R\ac\subseteq R$, which is the alcove in $R$ closest to the origin.
That is, for any $\alpha\in\Phi^+$ and $w_a\in\wa$ such that $w_aA_{\circ}\subseteq R$, we have $|k(w_R,\alpha)|\leq |k(w_a,\alpha)|$.
\end{mythrm}

\sk
We define $\wshi=\{w_R: R\text{ is a Shi region}\}$. The corresponding alcoves $w_R\ac$ we call \emph{Shi alcoves}. That is, we call an alcove a Shi alcove if it is the minimal alcove of the Shi region containing it.
\begin{mythrm}\label{thrm:alc} {\textnormal{\cite[Prop~7.3]{Shi1987}}}
The alcove $w_a\ac$ is a Shi alcove if and only if all floors of $w_a\ac$ are hyperplanes of the Shi arrangement.
\end{mythrm}

\sk
The following theorem is already known for dominant regions \cite[Prop~3.11]{Athan2005}.
\begin{mythrm}\label{thrm:floor} The floors of the minimal alcove $w_R\ac$ of a Shi region $R$ are exactly the floors of $R$.
\end{mythrm}

\sk
The following lemma describes what the Shi arrangement looks like in each chamber.
\begin{mylem}\label{thrm:cham} {\textnormal{\cite[Lemma~10.2]{ARR2012}}}
For $w\in W$, the hyperplanes of the Shi arrangement that intersect the chamber $wC$ are exactly those of the form $H_{w(\alpha)}^1$ where $\alpha\in\Phi^+$ and $w(\alpha)\in\Phi^+$.
\end{mylem}

\sk
Thus if $w_R\ac$ is a Shi alcove contained in the Weyl chamber $wC$, then by \reft{alc} and \refl{cham} all its floors are of the form $H_{w(\alpha)}^1$ where $\alpha\in\Phi^+$ and $w(\alpha)\in\Phi^+$.
So $w^{-1}w_R\ac$ is a dominant alcove and its floors are of the form $w^{-1}(H_{w(\alpha)}^1)=H_{\alpha}^1$ with $\alpha\in\Phi^+$. It is thus a Shi alcove by \reft{alc}.
Conversely, if $w_R\ac$ is a dominant Shi alcove and $w\in W$, then $ww_R\ac$ is a Shi alcove if and only if $w(\alpha)\in\Phi^+$ whenever $H_{\alpha}^1$ is a floor of $w_R\ac$.
Thus the map 
\begin{align*}
\Theta:w_R\mapsto (w^{-1}w_R,w)
\end{align*}
where $w_R\ac\subseteq wC$, is a bijection from $\wshi$ to the set of pairs $(w_R,w)$ such that $w_R\ac$ is a dominant Shi alcove, $w\in W$ and $w(\alpha)\in\Phi^+$ whenever $H_{\alpha}^1$ is a floor of $w_R\ac$.

\sk
Define a partial order on $\Phi^+$ by 
$\alpha\leq\beta$ if and only if $\beta-\alpha$ can be written as a linear combination of simple roots with nonnegative integer coefficients.
The set of positive roots $\Phi^+$ with this partial order is called the \emph{root poset}. 
It turns out that the map 
\begin{align*}
FL:R\mapsto\{\alpha\in\Phi^+: H_{\alpha}^1\text{ is a floor of }R\}
\end{align*}
is a bijection from the set of dominant Shi regions of $\Phi$ to the set of antichains in the root poset of $\Phi$. See \cite{Shi1997}.
Putting $R\mapsto w_R$, $\Theta$ and $FL$ together and using \reft{floor} we get that the map 
\begin{align*}
R\mapsto (A,w),
\end{align*}
where $R\subseteq wC$ and $A=w^{-1}(FL(R))$, is a bijection from the set of Shi regions to the set of pairs $(A,w)$ such that $A$ is an antichain in the root poset, $w\in W$ and $w(A)\subseteq\Phi^+$.
A similar bijection using ceilings instead of floors is given in \cite[Prop~10.3]{ARR2012}. 
\subsection{Types A and C}
If $\Phi$ is of type $A_{n-1}$, we take $V=\{(x_1,x_2,\ldots,x_n)\in\R^n:\sum_{i=1}^nx_i=0\}$, $\Phi=\{e_i-e_j: i\neq j\}$ and $\Phi^+=\{e_i-e_j: i<j\}$.
The Weyl group $W$ is the symmetric group $\S_n$ that acts on $V$ by permuting coordinates, $\Q=\{(x_1,x_2,\ldots,x_n)\in\Z^n:\sum_{i=1}^nx_i=0\}$, $r=n-1$ and $h=n$.

If $\Phi$ is of type $C_n$, we choose $V=\R^n$, $\Phi=\{e_i\pm e_j: i\neq j\}\cup\{\pm 2e_i: i\in[n]\}$ and $\Phi^+=\{e_i\pm e_j: i>j\}\cup\{2e_i: i\in[n]\}$.
The Weyl group $W$ is the hyperoctahedral group $\hypoct_n$ that acts on $V$ by permuting coordinates and changing signs, $\Q=\Z^n$, $r=n$ and $h=2n$.

%
%


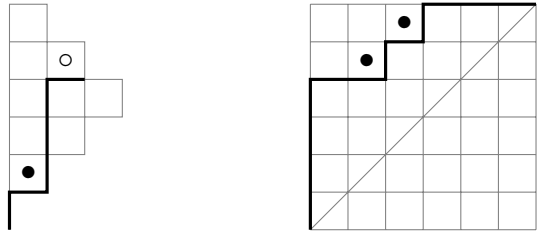
\begin{figure}[t]
\begin{center}
\begin{tikzpicture}[scale=.5]
\draw[gray] (0,1)--(1,1)--(1,6)--(0,6)--(0,0)
(0,2)--(2,2)--(2,5)--(0,5)
(0,3)--(3,3)--(3,4)--(0,4);
\draw[very thick] (0,0)--(0,1)--(1,1)--(1,4)--(2,4);
	\draw[xshift=5mm,yshift=5mm]
		(0,1) node{\large{$\bullet$}}
		(1,4) node{\large{$\circ$}};
\begin{scope}[xshift=8cm]
	\draw[gray] (0,0) grid (6,6)
		(0,0)--(6,6);
	\draw[very thick] (0,0)--(0,4)--(2,4)--(2,5)--(3,5)--(3,6)--(6,6);
	\draw (1.5,4.5) node{\large{$\bullet$}}
		(2.5,5.5) node{\large{$\bullet$}};
\end{scope}
\end{tikzpicture}
\caption{A ballot path $\beta\in\mathcal B_3$ (left) with one valley $(1,2)$ and two rises $2,3$, and a Dyck path $\pi\in\mathcal D_6$ (right) with valleys $(2,5),(3,6)$ and rises $1,2,3$. We have $A_{\beta}=\{e_3-e_2,2e_2\}$ and $A_{\pi}=\{e_2-e_5,e_3-e_6\}$}
\label{Figure:risevalley}
\end{center}
\end{figure}

\subsection{Lattice paths}

Denote by $\mathcal L_{m,n}$ the set of lattice paths from $(0,0)$ to $(m,n)$ consisting of $n$ North steps $N=(0,1)$ and $m$ East steps $E=(1,0)$. Let $\mathcal D_n$ denote the set of \emph{Dyck paths}, that is the subset of $\mathcal L_{n,n}$ consisting of the paths that never go below the main diagonal $x=y$. Let $\mathcal B_n$ denote the set of \emph{ballot paths}, that is the set of lattice paths starting at $(0,0)$, consisting of $2n$ North and/or East steps, and never going below the main diagonal.

A pattern of the form $NN$ is called \emph{rise}. A pattern $EN$ is called \emph{valley}. More precisely, let $\pi$ be any lattice path with steps $s_i\in\{N,E\}$. We say $i$ is a rise of $\pi$ if the $i$-th North step is followed by a North step. We say $(i,j)$ is a valley of $\pi$ if the $i$-th East step is followed by the $j$-th North step. See \reff{risevalley}.

If $\pi\in\D_n$ and $\Phi$ is of type $A_{n-1}$, define $A_{\pi}\subseteq\Phi^+$ by $e_i-e_j\in A_{\pi}$ if and only if $(i,j)$ is a valley of $\pi$.
Then the map $\pi\mapsto A_{\pi}$ is a bijection from $\D_n$ to the set of antichains in the root poset of $\Phi$.

If $\beta\in\B_n$ and $\Phi$ is of type $C_n$, define $A_{\beta}\subseteq\Phi^+$ by 
\begin{align*}
A_{\beta}&=\{e_i-e_j: i>j\text{ and }(n+1-i,n+1-j)\text{ is a valley of }\beta\}\\
&\qquad\cup\{e_i+e_j: i> j\text{ and }(n+1-i,j+n)\text{ is a valley of }\beta\}\\
&\qquad\cup\{2e_i:\text{the last step of $\beta$ is its $(n+1-i)$-th east step}\}.
\end{align*}
Then the map $\beta\mapsto A_{\beta}$ is a bijection from $\B_n$ to the set of antichains in the root poset of $\Phi$.
%
%

\section{Shi regions and parking functions}\label{sec:pf}

\subsection{Shi regions as diagonally labelled paths}\label{dlp}

A \emph{diagonally $\S_n$-labelled Dyck path} is a pair $(\pi,\sigma)$ of a Dyck path $\pi\in\mathcal D_n$ and a permutation $\sigma\in\S_n$ such that for each valley $(i,j)$ of $\pi$ we have $\sigma_i<\sigma_j$. See \reff{diagareaA}.
From the considerations at the end of Section \ref{Subsection:shi}, recall that regions of the Shi arrangement of type $A_{n-1}$ may be indexed by pairs $(A,\sigma)$ with $A$ an antichain in the root poset, $\sigma\in W=\S_n$ and $\sigma(A)\subseteq\Phi^+$.
\begin{myprop}\label{thrm:diagonalbijA} The map $(\pi,\sigma)\mapsto(A_{\pi},\sigma)$ is a bijection between diagonally labelled Dyck paths of length $n$ and regions of the Shi arrangement of type $A_{n-1}$.
\end{myprop}

We provide an interpretation of type $C$ Shi regions as diagonally labelled ballot paths. For any signed permutation $\sigma\in\hypoct_n$ we define $w^{\sigma}$ to be the word of length $2n$ given by $w^{\sigma}_i=\sigma(n+1-i)$ if $1\leq i\leq n$ and $w^{\sigma}_i=\sigma(n-i)$ if $n+1\leq i\leq2n$. For example if $n=3$ then $w^{\id}=321\bar1\bar2\bar3$.

A \emph{diagonally $\hypoct_n$-labelled ballot path} is a pair $(\beta,w^{\sigma})$ of a ballot path $\beta\in\mathcal B_n$ and a word $w^{\sigma}$ corresponding to a signed permutation $\sigma$ such that for each valley $(i,j)$ of $\beta$ we have $w^{\sigma}_i>w^{\sigma}_j$, and such that $0<w^{\sigma}_i$ if the final step of $\beta$ is its $i$-th East step. Hence, if we place the labels $w^{\sigma}$ in the diagonal then for each valley the label to its right will be smaller than the label below it. Moreover, if the path ends with an East step then the label below will be positive. See \reff{diagareaC}.


\begin{myprop}\label{thrm:diagonalbijC} The map $(\beta,w^\sigma)\mapsto(A_{\beta},\sigma)$ is a bijection between diagonally $\hypoct_n$-labelled ballot paths the regions of the Shi arrangement of type $C_n$.
\end{myprop}

\begin{figure}[t]
\begin{minipage}[t]{.45\linewidth}
\begin{center}
\begin{tikzpicture}[scale=.5]
\begin{scope}
	\fill[black!12] (1,2)--(2,2)--(2,3)--(3,3)--(3,5)--(2,5)--(2,4)--(1,4)--cycle;
	\draw[gray] (0,0) grid (5,5);
	\draw[very thick] (0,0)--(0,1)--(1,1)--(1,4)--(2,4)--(2,5)--(5,5);
	\draw[xshift=5mm,yshift=5mm]
		(0,0) node{\large{$1$}}
		(1,1) node{\large{$2$}}
		(2,2) node{\large{$3$}}
		(3,3) node{\large{$5$}}
		(4,4) node{\large{$4$}}
		(0,1) node{\large{$\bullet$}}
		(1,4) node{\large{$\bullet$}};
\end{scope}
\end{tikzpicture}
\caption{A diagonally labelled Dyck path $(\pi,\sigma)$, where $\sigma_1=1<\sigma_2=2$ and $\sigma_2=2<\sigma_5=4$. We have $\op{area}(\pi)=5$ and $\op{area'}(\pi,\sigma)=4$.}
\label{Figure:diagareaA}
\end{center}
\end{minipage}
\hfill
\begin{minipage}[t]{.45\linewidth}
\begin{center}
\begin{tikzpicture}[scale=.5]
\begin{scope}
	\fill[black!12] (1,2)--(2,2)--(2,3)--(1,3)--cycle;
	\draw[gray] (0,1)--(1,1)--(1,6)--(0,6)--(0,0)
		(0,2)--(2,2)--(2,5)--(0,5)
		(0,3)--(3,3)--(3,4)--(0,4);
	\draw[very thick] (0,0)--(0,2)--(1,2)--(1,4)--(2,4);
	\draw[xshift=5mm,yshift=5mm]
		(0,0) node{\large{$-1$}}
		(1,1) node{\large{$2$}}
		(2,2) node{\large{$-3$}}
		(3,3) node{\large{$3$}}
		(4,4) node{\large{$-2$}}
		(5,5) node{\large{$1$}}
		(1,4) node{\large{$\circ$}}
		(0,2) node{\large{$\bullet$}};
\end{scope}
\end{tikzpicture}
\caption{A diagonally labelled ballot path $(\beta,w^{\sigma})$, where $w_1=\sigma_{4-1}=-1>w_3=\sigma_{4-3}=-3$ and $0<w_2=\sigma_{4-2}=2$. We have $\op{area}(\beta)=4$ and $\op{area'}(\beta,w^{\sigma})=1$.}
\label{Figure:diagareaC}
\end{center}
\end{minipage}
\end{figure}

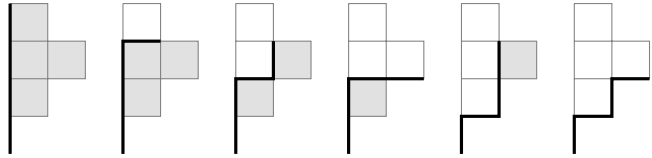
\begin{figure}[ht]
\begin{center}
\begin{tikzpicture}[scale=.5]
\fill[black!12] (0,1)--(1,1)--(1,2)--(2,2)--(2,3)--(1,3)--(1,4)--(0,4)--cycle;
\draw[gray] (0,0)--(0,4)
	(0,1)--(1,1)--(1,4)--(0,4)
	(0,2)--(2,2)--(2,3)--(0,3);
\draw[very thick] (0,0)--(0,4);
\begin{scope}[xshift=3cm]
\fill[black!12] (0,1)--(1,1)--(1,2)--(2,2)--(2,3)--(0,3)--(0,4)--cycle;
\draw[gray] (0,0)--(0,4)
	(0,1)--(1,1)--(1,4)--(0,4)
	(0,2)--(2,2)--(2,3)--(0,3);
\draw[very thick] (0,0)--(0,3)--(1,3);
\end{scope}
\begin{scope}[xshift=6cm]
\fill[black!12] (0,1)--(1,1)--(1,2)--(2,2)--(2,3)--(1,3)--(1,2)--(0,2)--cycle;
\draw[gray] (0,0)--(0,4)
	(0,1)--(1,1)--(1,4)--(0,4)
	(0,2)--(2,2)--(2,3)--(0,3);
\draw[very thick] (0,0)--(0,2)--(1,2)--(1,3);
\end{scope}
\begin{scope}[xshift=9cm]
\fill[black!12] (0,1)--(1,1)--(1,2)--(0,2)--cycle;
\draw[gray] (0,0)--(0,4)
	(0,1)--(1,1)--(1,4)--(0,4)
	(0,2)--(2,2)--(2,3)--(0,3);
\draw[very thick] (0,0)--(0,2)--(2,2);
\end{scope}
\begin{scope}[xshift=12cm]
\fill[black!12] (1,2)--(2,2)--(2,3)--(1,3)--cycle;
\draw[gray] (0,0)--(0,4)
	(0,1)--(1,1)--(1,4)--(0,4)
	(0,2)--(2,2)--(2,3)--(0,3);
\draw[very thick] (0,0)--(0,1)--(1,1)--(1,3);
\end{scope}
\begin{scope}[xshift=15cm]
\draw[gray] (0,0)--(0,4)
	(0,1)--(1,1)--(1,4)--(0,4)
	(0,2)--(2,2)--(2,3)--(0,3);
\draw[very thick] (0,0)--(0,1)--(1,1)--(1,2)--(2,2);
\end{scope}
\end{tikzpicture}
\caption{All ballot paths of length two and their area squares shaded gray.}
\label{Figure:areaC}
\end{center}
\end{figure}

\subsection{The area' statistic}

The area of a Dyck path is defined as the number of boxes strictly between the path and the main diagonal.
For example the Dyck path in \reff{diagareaA} has $\op{area}(\pi)=5$.
\cite{HagLoe2005} 
defined a related statistic $\op{area'}$ for diagonally labelled Dyck paths as follows. Consider a box strictly between the diagonal and the path $\pi$ in column $i$ and row $j$. This box contributes to $\op{area'}(\pi,\sigma)$ if and only if the label to its right is larger than the label below it, that is if and only if $\sigma_i<\sigma_j$. 
For example the labelled Dyck path in \reff{diagareaA} has $\op{area'}(\pi,\sigma)=4$ because the nonshaded box in the fifth row and fourth column does not contribute: $\sigma_4=5>\sigma_5=4$. 

\sk
The area of a ballot path is defined as the number of boxes ``below'' the path (see \reff{areaC}). We now define a type $C$ $\op{area'}$ statistic on diagonally labelled ballot paths. Let $(\beta,w^{\sigma})$ be such a path and consider a box below $\beta$ in column $i$ and row $j$. This box contributes to $\op{area'}(\beta,w^{\sigma}) $ if and only if the label to its right is smaller than the label below it, that is if and only if $w^{\sigma}_i>w^{\sigma}_j$. The labelled ballot path in \reff{diagareaC} has $\op{area'}(\beta,w^{\sigma})=1$ since the shaded box in the third row and second column is the only one contributing: $w^{\sigma}_2=2>w_3^{\sigma}=-3$. For example the box in the fourth row and second column does not contribute because $w^{\sigma}_2=2<w^{\sigma}_4=\sigma_{3-4}=3$.

\sk
Note that these statistics are the type $A$ and $C$ cases of the following uniform statistic.
Define the \emph{coheight} statistic on regions of the Shi arrangement of any irreducible root system $\Phi$ by
\begin{align*}
\op{coheight}(R)=|\Phi^+|-\#\text{ hyperplanes of the Shi arrangement separating }R\text{ from the origin}.
\end{align*}
Then the $\op{area'}$ statistics correspond to the $\op{coheight}$ statistic under the bijections in Section \ref{dlp}.

\subsection{Parking functions}
A vector $f=(f_1,\dots,f_n)$ with nonnegative integer entries is called a \emph{(classical) parking function} of length $n$ if there exists a permutation $\sigma\in\S_n$ such that $f_{\sigma(i)}\leq i-1$ for $1\leq i\leq n$. Equivalently, $f$ is a parking function if $\#\{j:f_j\leq i-1\}\geq i$ for all $1\leq i \leq n$.

There is a natural $\S_n$-isomorphism between the set of parking functions of length $n$ and the finite torus $\Q/(h+1)\Q$ of the root system of type $A_{n-1}$. Thus classical parking functions may be seen as objects of type $A$.

We define a \emph{type $C$ parking function} of length $n$ to be an integer vector $f=(f_1,\dots,f_n)$ where $-n\leq f_i\leq n$ for all $1\leq i\leq n$. Thus type $C$ parking functions of length $n$ are a natural set of representatives for the finite torus $\Q/(h+1)\Q=\Z^n/(2n+1)\Z^n$ of the root system of type $C_n$. 

\subsection{Vertically labelled paths}

Type $A$ parking functions are commonly represented as Dyck paths with labelled North steps \cite[Chap.~5]{Haglund2008}. An \emph{$\S_n$-labelled Dyck path} is a pair $(\pi,\sigma)$ of a Dyck path $\pi\in\mathcal D_n$ and a permutation $\sigma\in\S_n$ such that $\sigma_i<\sigma_{i+1}$ whenever $i$ is a rise of $\pi$. Thus, if the label $\sigma_i$ is placed in the box to the right of the $i$-th North step then labels increase along columns from bottom to top. For example in \reff{vertdinvA} we have $\sigma_1=1<\sigma_2=2<\sigma_3=4$.


\sk
We show how type $C$ parking functions can be regarded as labelled lattice paths in a similar fashion. An \emph{$\hypoct_n$-labelled path} $(\pi,\sigma)$ is a pair of a lattice path $\pi\in\mathcal L_{n,n}$ and a signed permutation $\sigma\in\mathfrak H_n$ such that $\sigma_i<\sigma_{i+1}$ whenever $i$ is a rise of $\pi$ and such that $0<\sigma_1$ if $\pi$ begins with a North step. Thus, if we place the label $\sigma_i$ to the left of the $i$-th North step then the labels increase along columns from bottom to top, and all labels in the zeroth column (that is left of the starting point) are positive. See \reff{verticalC}.

Given a parking function $f=(f_1,\dots,f_n)$ we obtain a labelled path as follows. For all $1\leq i\leq n$ if $f_i$ is non-negative, place the label $i$ in the $f_i$-th column. If $f_i$ is negative, place the label $-i$ in column $-f_i$. Rearrange the labels in each column in increasing order and draw a path as in \reff{verticalbijC}.

Conversely, let $(\pi,\sigma)$ be a labelled path. We define a parking function $g$ as follows. If a positive label $i$ occurs in the $j$-th column then set $g_i=j$. If a negative label $i$ occurs in the $j$-th column instead set $g_{-i}=-j$. In summary we have the following result.

\begin{myprop}\label{thrm:verticalbijC} The above correspondence defines a bijection between type $C$ parking functions of length $n$ and vertically $\hypoct_n$-labelled lattice paths.
\end{myprop}

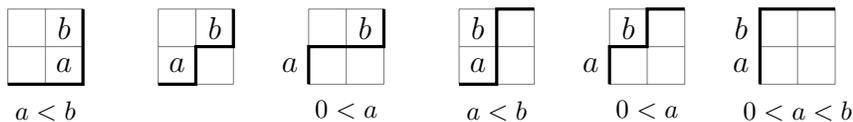
\begin{figure}[t]
\begin{center}
\begin{tikzpicture}[scale=.5]
\begin{scope}[xshift=0cm,yshift=0cm]
\draw[gray] (0,0) grid (2,2);
\draw[very thick] (0,0)--(2,0)--(2,2);
\draw (1.5,0.5) node{\large{$a$}};
\draw (1.5,1.5) node{\large{$b$}};
\draw (1,-.2) node[anchor=north]{{$a<b$}};
\end{scope}
\begin{scope}[xshift=4cm,yshift=0cm]
\draw[gray] (0,0) grid (2,2);
\draw[very thick] (0,0)--(1,0)--(1,1)--(2,1)--(2,2);
\draw (0.5,0.5) node{\large{$a$}};
\draw (1.5,1.5) node{\large{$b$}};
\end{scope}
\begin{scope}[xshift=8cm,yshift=0cm]
\draw[gray] (0,0) grid (2,2);
\draw[very thick] (0,0)--(0,1)--(2,1)--(2,2);
\draw (-.5,0.5) node{\large{$a$}};
\draw (1.5,1.5) node{\large{$b$}};
\draw (1,-.2) node[anchor=north]{{$0<a$}};
\end{scope}
\begin{scope}[xshift=12cm]
\begin{scope}[xshift=0cm]
\draw[gray] (0,0) grid (2,2);
\draw[very thick] (0,0)--(1,0)--(1,2)--(2,2);
\draw (0.5,.5) node{\large{$a$}};
\draw (0.5,1.5) node{\large{$b$}};
\draw (1,-.2) node[anchor=north]{{$a<b$}};
\end{scope}
\begin{scope}[xshift=4cm]
\draw[gray] (0,0) grid (2,2);
\draw[very thick] (0,0)--(0,1)--(1,1)--(1,2)--(2,2);
\draw (-.5,.5) node{\large{$a$}};
\draw (0.5,1.5) node{\large{$b$}};
\draw (1,-.2) node[anchor=north]{{$0<a$}};
\end{scope}
\begin{scope}[xshift=8cm]
\draw[gray] (0,0) grid (2,2);
\draw[very thick] (0,0)--(0,2)--(2,2);
\draw (-.5,.5) node{\large{$a$}};
\draw (-.5,1.5) node{\large{$b$}};
\draw (1,-.2) node[anchor=north]{{$0<a<b$}};
\end{scope}
\end{scope}
\end{tikzpicture}
\caption{All six paths in $\mathcal L_{2,2}$ and the conditions on their labellings.}
\label{Figure:verticalC}
\end{center}
\end{figure}

\begin{figure}[t]
\begin{center}
\begin{tikzpicture}[scale=.5]
\begin{scope}
	\draw[gray] (0,0) grid (4,4);
	\draw[xshift=5mm,yshift=5mm]
		(-1,0) node{\large{$2$}}
		(0,0) node{\large{$-3$}}
		(3,0) node{\large{$1$}}
		(3,1) node{\large{$-4$}};
	\draw[->,thick] (5,2)--(5.8,2);
\end{scope}
\begin{scope}[xshift=7cm]
	\draw[gray] (0,0) grid (4,4);
	\draw[very thick] (0,0)--(0,1)
		(1,1)--(1,2)
		(4,2)--(4,4);
	\draw[xshift=5mm,yshift=5mm]
		(-1,0) node{\large{$2$}}
		(0,1) node{\large{$-3$}}
		(3,2) node{\large{$-4$}}
		(3,3) node{\large{$1$}};
	\draw[->,thick] (5,2)--(5.8,2);
\end{scope}
\begin{scope}[xshift=14cm]
	\draw[gray] (0,0) grid (4,4);
	\draw[very thick] (0,0)--(0,1)--(1,1)--(1,2)--(4,2)--(4,4);
	\draw (-.5,.5) node{\large{$2$}}
		(0.5,1.5) node{\large{$-3$}}
		(3.5,2.5) node{\large{$-4$}}
		(3.5,3.5) node{\large{$1$}};
\end{scope}
\end{tikzpicture}
\caption{Constructing an $\hypoct_4$-labelled path from the parking function $f=(4,0,-1,-4)$.}
\label{Figure:verticalbijC}
\end{center}
\end{figure}
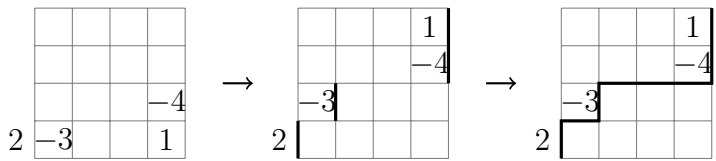

\begin{figure}[t]
\begin{minipage}[t]{.45\linewidth}
\begin{center}
\begin{tikzpicture}[scale=.5]
\begin{scope}
	\draw[gray] (0,0) grid (5,5);
	\draw[very thick] (0,0)--(0,3)--(2,3)--(2,4)--(3,4)--(3,5)--(5,5);
	\draw[xshift=5mm,yshift=5mm]
		(0,0) node{\large{1}}
		(0,1) node{\large{2}}
		(0,2) node{\large{4}}
		(2,3) node{\large{3}}
		(3,4) node{\large{5}};
\end{scope}
\end{tikzpicture}
\caption{A vertically labelled Dyck path $(\pi,\sigma)$ with area vector $(0,1,2,1,1)$ and $\op{dinv}(\pi)=5$ and $\op{dinv'}(\pi,\sigma)=4$.}
\label{Figure:vertdinvA}
\end{center}
\end{minipage}
\hfill
\begin{minipage}[t]{.45\linewidth}
\begin{center}
\begin{tikzpicture}[scale=.5]
\begin{scope}
	\draw[gray] (0,0) grid (6,6);
	\draw[very thick] (0,0)--(0,1)--(4,1)--(4,6)--(6,6);
	\draw[xshift=5mm,yshift=5mm]
		(-1,0) node{\large{$1$}}
		(3,1) node{\large{$-5$}}
		(3,2) node{\large{$-4$}}
		(3,3) node{\large{$2$}}
		(3,4) node{\large{$3$}}
		(3,5) node{\large{$6$}};
\end{scope} 
\end{tikzpicture}
\caption{A vertically labelled path $(\beta,\sigma)$ with area vector $(1,-2,-1,0,1,2)$ and $\op{dinv}(\beta)=9$ and $\op{dinv'}(\pi,\sigma)=6$.}
\label{Figure:vertdinvC}
\end{center}
\end{minipage}
\end{figure}

\subsection{The dinv' statistic}

The $\op{dinv}$ statistic was first defined by Haiman to provide an (at the time conjectural) combinatorial model for the $q,t$-Catalan numbers \cite[Chap. 3]{Haglund2008}. For each Dyck path $\pi\in\mathcal D_n$ define the area vector $(a_1,a_2,\dots,a_n)$ by letting $a_i$ be the number of boxes in the $i$-th row, strictly between $\pi$ and the main diagonal. For example the Dyck path in \reff{vertdinvA} has area vector $(0,1,2,1,1)$. The $\op{dinv}$ statistic is defined as
\begin{align*}
\op{dinv}(\pi)=\#\big\{(i,j):i<j,a_i=a_j\big\}
+\#\big\{(i,j):i<j,a_i=a_j+1\big\}.
\end{align*}
A pair $(i,j)$ contributing to $\op{dinv}$ is called a \emph{diagonal inversion}. 
\cite{HagLoe2005} defined a generalised statistic $\op{dinv'}$ on vertically $\S_n$-labelled Dyck paths. A pair $(i,j)$ with $a_i=a_j$ contributes if and only if $\sigma_i<\sigma_j$. On the other hand, a pair $(i,j)$ with $a_i=a_j+1$ contributes if and only if $\sigma_i>\sigma_j$. Compare with \reff{vertdinvA}.

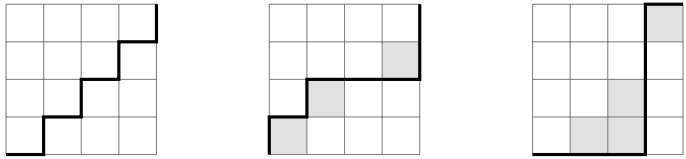
\begin{figure}[t]
\begin{center}
\begin{tikzpicture}[scale=.5]
\draw[gray] (0,0) grid (4,4);
\draw[very thick] (0,0)--(1,0)--(1,1)--(2,1)--(2,2)--(3,2)--(3,3)--(4,3)--(4,4);
\begin{scope}[xshift=7cm]
\fill[black!12] (0,0)--(1,0)--(1,2)--(2,2)--(2,1)--(0,1)--cycle
	(3,2)--(4,2)--(4,3)--(3,3)--cycle;
\draw[gray] (0,0) grid (4,4);
\draw[very thick] (0,0)--(0,1)--(1,1)--(1,2)--(4,2)--(4,4);
\end{scope}
\begin{scope}[xshift=14cm]
\fill[black!12] (1,0)--(3,0)--(3,2)--(2,2)--(2,1)--(1,1)--cycle
	(3,3)--(4,3)--(4,4)--(3,4)--cycle;
\draw[gray] (0,0) grid (4,4);
\draw[very thick] (0,0)--(3,0)--(3,4)--(4,4);
\end{scope}
\end{tikzpicture}
\caption{Paths with type $C$ area vectors $(0,0,0,0),(1,1,-1,0)$ and $(-2,-1,0,1)$.}
\label{Figure:areaseqC}
\end{center}
\end{figure}

\sk
We define an area vector and a $\op{dinv}$ statistic of type $C$ for lattice paths $\pi\in\mathcal L_{n,n}$. The area vector $(a_1,a_2,\dots,a_n)$ is given by $a_i=i-b_i$ where $b_i$ is the number of boxes in the $i$-th row left of $\pi$. See Figures \ref{Figure:vertdinvC} and \ref{Figure:areaseqC}. Moreover, we define
\begin{align*}
\op{dinv}(\pi)
&= \#\big\{(i,j):i<j,a_i=a_j\big\} + \#\big\{(i,j):i<j,a_i=a_j+1\big\} \\
&\qquad+ \#\big\{(i,j):i<j,a_i=-a_j\big\} + \#\big\{(i,j):i<j,a_i=-a_j+1\big\} \\
&\qquad+ \#\big\{i:a_i=0\big\}.
\end{align*}

Next, we define a $\op{dinv'}$ statistic for vertically $\hypoct_n$-labelled lattice path $(\pi,\sigma)$. As in type $A$ above, a pair $(i,j)$ of candidate rows contributes if and only if the labels $\sigma_i$ and $\sigma_j$ satisfy a certain inequality. More precisely,
\begin{align*}
\notag
\op{dinv'}(\pi,\sigma)
&= \#\big\{(i,j):i<j,a_i=a_j,\sigma_i<\sigma_j\big\}
+ \#\big\{(i,j):i<j,a_i=a_j+1,\sigma_i>\sigma_j\big\} \\
\notag
&\qquad+ \#\big\{(i,j):i<j,a_i=-a_j,\sigma_i<-\sigma_j\big\}
+ \#\big\{(i,j):i<j,a_i=-a_j+1,\sigma_i>-\sigma_j\big\}\\ 
&\qquad+ \#\big\{i:a_i=0,\sigma_i<0\big\}.
\end{align*}
Note that the two definitions of $\op{dinv}(\pi)$ agree if $\pi$ is a Dyck path.
\sk
Consider the $\hypoct_6$-labelled path $(\pi,\sigma)$ in \reff{vertdinvC}. Its area vector is given by $(1,-2,-1,0,1,2)$. There is one diagonal inversion of type $a_i=a_j$, namely $(1,5)$, one diagonal inversion of type $a_i=a_j+1$, namely $(1,4)$, three diagonal inversions of type $a_i=-a_j$, namely $(1,3),(2,6)$ and $(3,5)$, three diagonal inversions of type $a_i=-a_j+1$, namely $(1,4),(3,6)$ and $(4,5)$, and one row of length zero, namely $i=4$. In total we have 9 diagonal inversions, so $\op{dinv}(\pi)=9$. Note that the inversion $(1,4)$ is counted twice!

If we wish to take labels into account we find that $\sigma_4=2>0$, so the row of length zero does not contribute. Moreover $\sigma_1=1<\sigma_4=2$, thus $(1,4)$ is not a d'-inversion of type $a_i=a_j+1$, and $\sigma_2=-5>-\sigma_6=-6$, thus $(2,6)$ is not a d'-inversion of type $a_i=-a_j$. The labels of all other diagonal inversions fit our requirements, so $\op{dinv'}(\pi,\sigma)=6$.

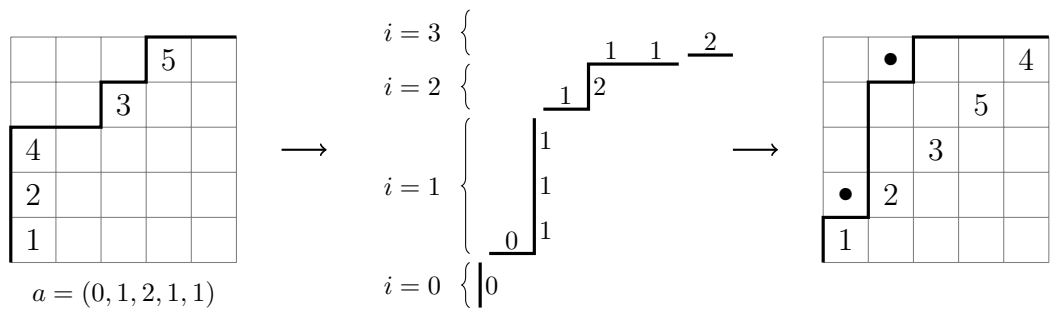
\begin{figure}[t]
\begin{center}
\begin{tikzpicture}[scale=.6]
\begin{scope}
	\draw[gray] (0,0) grid (5,5);
	\draw[very thick] (0,0)--(0,3)--(2,3)--(2,4)--(3,4)--(3,5)--(5,5);
	\draw (2.5,-.2) node[anchor=north]{$a=(0,1,2,1,1)$};
	\draw[->,thick] (6,2.5)--(7,2.5);
	\draw[xshift=5mm,yshift=5mm]
		(0,0) node{\large{$1$}}
		(0,1) node{\large{$2$}}
		(0,2) node{\large{$4$}}
		(2,3) node{\large{$3$}}
		(3,4) node{\large{$5$}};
\end{scope}
\begin{scope}[xshift=8.6cm,yshift=-1cm]
	\begin{scope}[xshift=1.8cm,yshift=0cm]
		\draw[very thick] (0,0)--(0,1);
		\draw[xshift=-1mm,yshift=5mm,anchor=west]
			(0,0) node{$0$};
		\draw[decorate,decoration={brace,amplitude=3pt,raise=4pt}] (0,0)--(0,1);
		\draw (-1.5,0.5) node{$i=0$};

	\end{scope}
	\begin{scope}[xshift=2cm,yshift=1.2cm]
		\draw[very thick] (0,0)--(1,0)--(1,3);
		\draw[xshift=5mm,yshift=-1mm,anchor=south]
			(0,0) node{$0$};
		\draw[xshift=-1mm,yshift=5mm,anchor=west]
			(1,0) node{$1$}
			(1,1) node{$1$}
			(1,2) node{$1$};
		\draw[decorate,decoration={brace,amplitude=3pt,raise=4pt}] (-.2,0)--(-.2,3);
		\draw (-1.7,1.5) node{$i=1$};
	\end{scope}
	\begin{scope}[xshift=3.2cm,yshift=4.4cm]
		\draw[very thick] (0,0)--(1,0)--(1,1)--(3,1);
		\draw[xshift=5mm,yshift=-1mm,anchor=south]
			(0,0) node{$1$}
			(1,1) node{$1$}
			(2,1) node{$1$};
		\draw[xshift=-1mm,yshift=5mm,anchor=west]
			(1,0) node{$2$};
		\draw[decorate,decoration={brace,amplitude=3pt,raise=4pt}] (-1.4,0)--(-1.4,1);
		\draw (-2.9,.5) node{$i=2$};
	\end{scope}
	\begin{scope}[xshift=6.4cm,yshift=5.6cm]
		\draw[very thick] (0,0)--(1,0);
		\draw[xshift=5mm,yshift=-1mm,anchor=south]
			(0,0) node{$2$};
		\draw[decorate,decoration={brace,amplitude=3pt,raise=4pt}] (-4.6,0)--(-4.6,1);
		\draw (-6.1,.5) node{$i=3$};
	\end{scope}
\end{scope}
\begin{scope}[xshift=18cm]
	\draw[gray] (0,0) grid (5,5);
	\draw[very thick] (0,0)--(0,1)--(1,1)--(1,4)--(2,4)--(2,5)--(5,5);
	\draw[->,thick] (-2,2.5)--(-1,2.5);
	\draw[xshift=5mm,yshift=5mm]
		(0,0) node{\large{$1$}}
		(1,1) node{\large{$2$}}
		(2,2) node{\large{$3$}}
		(3,3) node{\large{$5$}}
		(4,4) node{\large{$4$}}
		(0,1) node{\large{$\bullet$}}
		(1,4) node{\large{$\bullet$}};
\end{scope}
\end{tikzpicture}
\caption{The classical zeta map: A vertically labelled Dyck path $(\pi,\sigma)$ (left), the construction of $\zeta(\pi)$ (middle), and $\zeta(\pi,\sigma)$ (right).}
\label{Figure:zetaA}
\end{center}
\end{figure}

\section{The zeta map}\label{sec:zeta}

The original zeta map is a bijection $\zeta:\mathcal D_n\to\mathcal D_n$ on Dyck paths and appears in a paper of 
\cite{AKOP2002}. A more explicit treatment including the compatibility with the statistics on Dyck paths defined in the previous sections can be found in \cite[Thrm.~3.15]{Haglund2008}. Let us start by recalling the definition of the zeta map.

Given a Dyck path $\pi\in\mathcal D_n$ with area vector $(a_1,a_2,\dots,a_n)$, set $i=0$ and place your pen at $(0,0)$. Now read the area vector from left to right drawing an East step for each $i-1$ you encounter and a North step for each $i$. Replace $i$ by $i+1$ and repeat until you reach the point $(n,n)$. See \reff{zetaA}.

\sk
We describe a bijection $\zeta_C:\mathcal L_{n,n}\to\mathcal{B}_n$ which is an analogue of the classical zeta map.

Given a path $\pi\in\mathcal L_{n,n}$ with type $C$ area vector $(a_1,a_2,\ldots,a_n)$,
set $i=n$ and start with your pen at $(0,0)$.
Read the area vector from left to right drawing an East step for each $-i-1$ you encounter and a North step for each $-i$. Then read the area vector from right to left drawing an East step for each $i+1$ you encounter and a North step for each $i$. Now replace $i$ by $i-1$ and repeat the process until $2n$ steps are drawn. See \reff{zetaC}.


\begin{figure}[t]
\begin{center}
\begin{tikzpicture}[scale=.55]
\begin{scope}[yshift=1cm]
	\draw[gray] (0,0) grid (6,6);
	\draw[very thick] (0,0)--(0,1)--(4,1)--(4,6)--(6,6);
	\draw[xshift=5mm,yshift=5mm]
		(-1,0) node{$1$}
		(3,1) node{$-5$}
		(3,2) node{$-4$}
		(3,3) node{$2$}
		(3,4) node{$3$}
		(3,5) node{$6$};
	\draw (3,-.2) node[anchor=north]{$a=(1,-2,-1,0,1,2)$};
	\draw[->,thick] (6.8,2.5)--(7.6,2.5);
\end{scope}
\begin{scope}[xshift=1cm]
	\begin{scope}[xshift=9.8cm,yshift=0cm]
		\draw[very thick] (0,0)--(0,2);
		\draw[xshift=-1mm,yshift=5mm,anchor=west]
			(0,0) node{$-2$}
			(0,1) node{$2$};
		\draw[decorate,decoration={brace,amplitude=3pt,raise=4pt}] (0,0)--(0,2);
		\draw (-1.5,1) node{$i=2$};
	\end{scope}
	\begin{scope}[xshift=10cm,yshift=2.2cm]
		\draw[very thick] (0,0)--(1,0)--(1,1)--(2,1)--(2,3);
		\draw[xshift=5mm,yshift=-1mm,anchor=south]
			(0,0) node{$-2$}
			(1,1) node{$2$};
		\draw[xshift=-1mm,yshift=5mm,anchor=west]
			(1,0) node{$-1$}
			(2,1) node{$1$}
			(2,2) node{$1$};
		\draw[decorate,decoration={brace,amplitude=3pt,raise=4pt}] (-.2,0)--(-.2,3);
		\draw (-1.7,1.5) node{$i=1$};
	\end{scope}
	\begin{scope}[xshift=12.2cm,yshift=5.4cm]
		\draw[very thick] (0,0)--(1,0)--(1,1)--(2,1)--(2,2)--(3,2);
		\draw[xshift=5mm,yshift=-1mm,anchor=south]
			(0,0) node{$-1$}
			(1,1) node{$1$}
			(2,2) node{$1$};
		\draw[xshift=-1mm,yshift=5mm,anchor=west]
			(1,0) node{$0$}
			(2,1) node{$0$};
		\draw[decorate,decoration={brace,amplitude=3pt,raise=4pt}] (-2.4,0)--(-2.4,2);
		\draw (-3.9,1) node{$i=0$};
	\end{scope}
\end{scope}
\begin{scope}[xshift=18cm,yshift=-.5cm]
	\fill[black!12] (1,2)--(2,2)--(2,3)--(3,3)--(3,4)--(4,4)--(4,5)--(5,5)--(5,6)--(3,6)--(3,5)--(2,5)--(2,3)--(1,3)--cycle;
	\draw[gray]
		(0,1)--(1,1)--(1,8)--(0,8)--(0,0)
		(0,2)--(2,2)--(2,8)--(0,8)
		(0,3)--(3,3)--(3,8)--(0,8)
		(0,4)--(4,4)--(4,8)--(0,8)
		(0,5)--(5,5)--(5,8)--(0,8)
		(0,6)--(6,6)--(6,7)--(0,7);
	\draw[very thick]
		(0,0)--(0,2)--(1,2)--(1,3)--(2,3)--(2,5)--(3,5)--(3,6)--(4,6)--(4,7)--(5,7);
	\draw[->,thick] (-1.6,4)--(-.8,4);
	\draw[xshift=5mm,yshift=5mm]
		(0,0) node{$5$}
		(1,1) node{$6$}
		(2,2) node{$4$}
		(3,3) node{$3$}
		(4,4) node{$1$}
		(5,5) node{$-2$}
		(6,6) node{$2$}
		(0,2) node{\large{$\bullet$}}
		(1,3) node{\large{$\bullet$}}
		(2,5) node{\large{$\bullet$}}
		(3,6) node{\large{$\bullet$}}
		(4,7) node{\large{$\circ$}};
	\draw[thick,loosely dotted] (7.1,7.1)--(7.7,7.7);
\end{scope}
\end{tikzpicture}
\caption{The type $C$ zeta map: A vertically labelled lattice path $(\pi,\sigma)$ (left), the construction of $\zeta_C(\pi)$ (middle), and $\zeta(\pi,\sigma)$ (right). Note that $\op{dinv'}(\pi,\sigma)=6=\op{area'}(\zeta(\pi,\sigma))$.
}
\label{Figure:zetaC}
\end{center}
\end{figure}
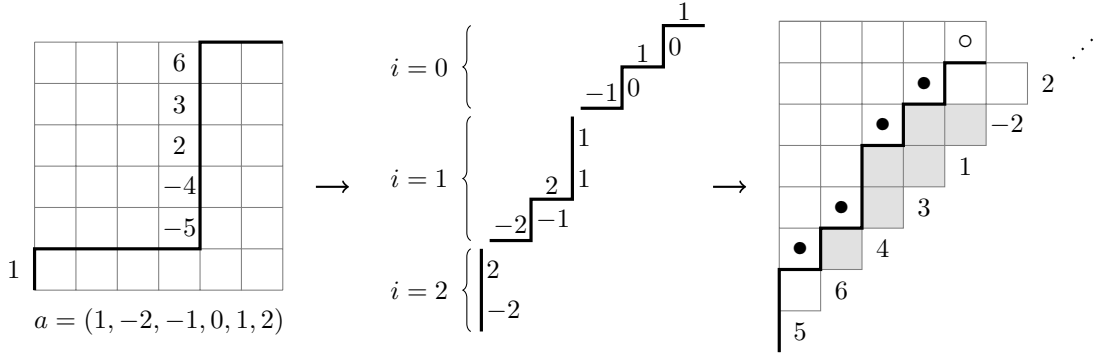

It is clear from the construction that $\zeta_C(\pi)$ 
never goes below the main diagonal. Moreover when $\pi$ is a Dyck path, then $\zeta_C(\pi)$ is just the reverse path of $\zeta(\pi)$. In particular $\zeta_C$ sends Dyck paths to Dyck paths.

The following is our first main result.

\begin{mythrm}\label{thrm:zetaC} The map $\zeta_C:\mathcal L_{n,n}\to\mathcal B_n$ is a bijection such that $\op{dinv}_C(\pi)=\op{area}(\zeta_C(\pi))$.
\end{mythrm}

The zeta map can be inverted using the bounce path of a ballot path. A detailed proof will appear in the full version.

\subsection{The Haglund--Loehr zeta map}

\cite{HagLoe2005} extended the classical zeta map to a bijection from vertically labelled Dyck paths to diagonally labelled Dyck paths that sends the $\op{dinv'}$ statistic to the $\op{area'}$ statistic. We start out by recalling their definition. If $(\pi,\sigma)$ is a vertically $\S_n$-labelled path, then $\zeta(\pi,\sigma)$ is simply the diagonal labelling of $\zeta(\pi)$ obtained as follows. For $i=0,1,\dots,n$ read the labels of rows with area equal to $i$ from bottom to top and insert them in the diagonal. Compare with \reff{zetaA}.

\sk
Similarly, in type $C$ we start with a vertically $\hypoct_n$-labelled path $(\pi,\sigma)$ and construct a diagonally labelled ballot path $\zeta(\pi,\sigma)=(\beta,w)$. The ballot path is given by $\beta=\zeta_C(\pi)$. The labelling is obtained as follows. For $i=n,n-1,\dots,1$ read the labels of the rows with area $i$ from top to bottom and insert them in the diagonal, then read the labels of rows with area equal to $-i+1$ from bottom to top and insert their negatives in the diagonal. In the end complement the $n$ labels by adding their negatives in reverse order. See Figures \ref{Figure:zetaC}, \ref{Figure:HLzetaC} and \ref{Figure:HLzetaC2}.

\begin{figure}[t]
\begin{minipage}[t]{.46\linewidth}
\begin{center}
\begin{tikzpicture}[scale=.5]
\begin{scope}
	\draw[gray] (0,1) grid (3,4);
	\draw[very thick] (0,1)--(0,2)--(2,2)--(2,4)--(3,4);
	\draw[xshift=5mm,yshift=5mm]
		(-1,1) node {\large{$2$}}
		(1,2) node {\large{$-1$}}
		(1,3) node {\large{$3$}};
	\draw[thick,->] (3.6,2)--(4.4,2);
	\draw
		(4,2.2) node[anchor=south]{\large{$\zeta$}}
		(1.5,.8) node[anchor=north]{{$(\pi,\sigma)$}};
\end{scope}
\begin{scope}[xshift=5cm]
	\fill[black!12] (0,1)--(1,1)--(1,2)--(2,2)--(2,3)--(3,3)--(3,4)--(1,4)--(1,3)--(0,3)--cycle;
	\draw[gray] (0,1)--(1,1)--(1,6)--(0,6)--(0,0)
		(0,2)--(2,2)--(2,5)--(0,5)
		(0,3)--(3,3)--(3,4)--(0,4);
	\draw[very thick] (0,0)--(0,3)--(1,3)--(1,4)--(2,4);
	\draw[xshift=5mm,yshift=5mm]
		(0,0) node{\large{$3$}}
		(1,1) node{\large{$2$}}
		(2,2) node{\large{$1$}}
		(3,3) node{\large{$-1$}}
		(4,4) node{\large{$-2$}}
		(5,5) node{\large{$-3$}}
		(1,4) node{\large{$\circ$}}
		(0,3) node{\large{$\bullet$}};
	\draw (3,.8) node[anchor=north]{{$(\beta,w)$}};
\end{scope}
\end{tikzpicture}
\caption{The area vector of $\pi$ is $(1,0,1)$. The diagonal inversions are $(1,2),(1,2),(1,3),(2,3),i=2$, and $\op{dinv'}(\pi,\sigma)=5=\op{area'}(\beta,w)$.}
\label{Figure:HLzetaC}
\end{center}
\end{minipage}
\hfill
\begin{minipage}[t]{.48\linewidth}
\begin{center}
\begin{tikzpicture}[scale=.5]
\begin{scope}
	\draw[gray] (0,1) grid (3,4);
	\draw[very thick] (0,1)--(0,2)--(3,2)--(3,4);
	\draw[xshift=5mm,yshift=5mm]
		(-1,1) node {\large{$2$}}
		(2,2) node {\large{$1$}}
		(2,3) node {\large{$3$}};
	\draw (1.5,.8) node[anchor=north]{{$(\pi,\sigma)$}};
	\draw[thick,->] (3.6,2)--(4.4,2);
	\draw (4,2.2) node[anchor=south]{\large{$\zeta$}};
\end{scope}
\begin{scope}[xshift=5cm]
	\fill[black!12] (1,2)--(2,2)--(2,3)--(1,3)--cycle;
	\draw[gray] (0,1)--(1,1)--(1,6)--(0,6)--(0,0)
		(0,2)--(2,2)--(2,5)--(0,5)
		(0,3)--(3,3)--(3,4)--(0,4);
	\draw[very thick] (0,0)--(0,2)--(1,2)--(1,4)--(2,4);
	\draw[xshift=5mm,yshift=5mm]
		(0,0) node{\large{$-1$}}
		(1,1) node{\large{$2$}}
		(2,2) node{\large{$-3$}}
		(3,3) node{\large{$3$}}
		(4,4) node{\large{$-2$}}
		(5,5) node{\large{$1$}}
		(1,4) node{\large{$\circ$}}
		(0,2) node{\large{$\bullet$}};
	\draw (3,.8) node[anchor=north]{{$(\beta,w)$}};
\end{scope}
\end{tikzpicture}
\caption{The area vector of $\pi$ is $(1,-1,0)$. The diagonal inversions are $(1,2),(1,3),(1,3),i=3$ but only $(1,3)$ (as inversion of type $a_i=-a_j+1$) contributes to $\op{dinv'}(\pi,\sigma)=1$.}
\label{Figure:HLzetaC2}
\end{center}
\end{minipage}
\end{figure}

\sk
The theorem below is the main result of this paper.

\begin{mythrm}\label{thrm:HLzetaC} The type $C$ zeta map is a bijection from vertically $\hypoct_n$-labelled paths to diagonally $\hypoct_n$-labelled ballot paths that sends the $\op{dinv'}$ statistic to the $\op{area'}$ statistic.
\end{mythrm}

Combining \reft{HLzetaC} with Propositions \ref{thrm:diagonalbijC} and \ref{thrm:verticalbijC} we obtain a new proof of the well known result that the Shi arrangement of type $C_n$ has $(2n+1)^n$ regions.

\subsection{The zeta maps via valleys}

The Haglund--Loehr zeta map has another simple description given in \cite{ALW2014}.
Let us fix the following convention. If $(\pi,\sigma)$ is a vertically $\S_n$-labelled Dyck path with rise $i$ then we say the rise has label $(\sigma_i,\sigma_{i+1})$. If $(\pi,\sigma)$ is a diagonally labelled Dyck path with valley $(i,j)$ then we say the valley has label $(\sigma_i,\sigma_j)$.

The image $\zeta(\pi,\sigma)$ of a vertically labelled Dyck path under the zeta map can now be defined as follows. First insert the diagonal labelling as described in the previous section. The Dyck path is the unique path which has a valley labelled $(a,b)$ if and only if $(\pi,\sigma)$ has a rise labelled $(a,b)$. 
See \reff{zetaA}.

\sk
There is a similar description of the zeta map in type $C$ provided by the following proposition. If $(\pi,\sigma)$ is a vertically $\hypoct_n$-labelled path with rise $i$ then we say the rise has label $(\sigma_i,\sigma_{i+1})$. If $(\beta,w)$ is a diagonally labelled ballot path with valley $(i,j)$ then we say the valley has label $(w_i,w_j)$.

\begin{myprop}\label{thrm:zetaCvalleys} Let $(\pi,\sigma)$ be a vertically $\hypoct_n$-labelled lattice path. Then $\zeta(\pi,\sigma)$ has valley labelled $(a,b)$ if and only if $(\pi,\sigma)$ has a rise labelled $(b,a)$ or $(-a,-b)$. Moreover, $\zeta(\pi,\sigma)$ ends with an East step in the same column as label $a$ if and only if $(\pi,\sigma)$ begins with a North step labelled $a$.
\end{myprop}

Compare with Figures \ref{Figure:zetaC}, \ref{Figure:HLzetaC} and \ref{Figure:HLzetaC2}.

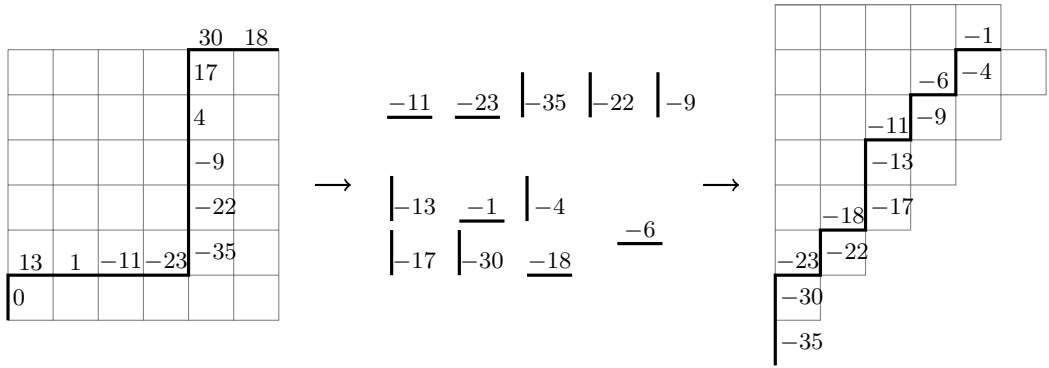
\begin{figure}[ht]
\begin{center}
\begin{tikzpicture}[scale=.6]
\begin{scope}
	\draw[gray] (0,0) grid (6,6);
	\draw[very thick] (0,0)--(0,1)--(4,1)--(4,6)--(6,6);
	\draw[xshift=5mm,yshift=-1mm,anchor=south]
		(0,1) node{\small{$13$}}
		(1,1) node{\small{$1$}}
		(2,1) node{\small{$-11$}}
		(3,1) node{\small{$-23$}}
		(4,6) node{\small{$30$}}
		(5,6) node{\small{$18$}};
	\draw[xshift=-1mm,yshift=5mm,anchor=west]
		(0,0) node{\small{$0$}}
		(4,1) node{\small{$-35$}}
		(4,2) node{\small{$-22$}}
		(4,3) node{\small{$-9$}}
		(4,4) node{\small{$4$}}
		(4,5) node{\small{$17$}};
	\draw[->,thick] (6.8,3)--(7.6,3);
\end{scope}
\begin{scope}[xshift=8.4cm,yshift=45mm]
	\begin{scope}
	\draw[very thick](0,0)--(1,0);
	\draw[xshift=5mm,yshift=-1mm,anchor=south] (0,0) node{\small{$-11$}};
	\end{scope}
	\begin{scope}[xshift=15mm]
	\draw[very thick](0,0)--(1,0);
	\draw[xshift=5mm,yshift=-1mm,anchor=south] (0,0) node{\small{$-23$}};
	\end{scope}
	\begin{scope}[xshift=30mm]
	\draw[very thick](0,0)--(0,1);
	\draw[xshift=5mm,yshift=-1mm,anchor=south] (0,0) node{\small{$-35$}};
	\end{scope}
	\begin{scope}[xshift=45mm]
	\draw[very thick](0,0)--(0,1);
	\draw[xshift=5mm,yshift=-1mm,anchor=south] (0,0) node{\small{$-22$}};
	\end{scope}
	\begin{scope}[xshift=60mm]
	\draw[very thick](0,0)--(0,1);
	\draw[xshift=5mm,yshift=-1mm,anchor=south] (0,0) node{\small{$-9$}};
	\end{scope}
\end{scope}
\begin{scope}[xshift=85mm,yshift=1cm]
	\begin{scope}[yshift=12mm]
	\draw[very thick](0,0)--(0,1);
	\draw[xshift=5mm,yshift=-1mm,anchor=south] (0,0) node{\small{$-13$}};
	\end{scope}
	\begin{scope}[xshift=15mm,yshift=12mm]
	\draw[very thick](0,0)--(1,0);
	\draw[xshift=5mm,yshift=-1mm,anchor=south] (0,0) node{\small{$-1$}};
	\end{scope}
	\begin{scope}[xshift=30mm,yshift=12mm]
	\draw[very thick](0,0)--(0,1);
	\draw[xshift=5mm,yshift=-1mm,anchor=south] (0,0) node{\small{$-4$}};
	\end{scope}
	\begin{scope}
	\draw[very thick](0,0)--(0,1);
	\draw[xshift=5mm,yshift=-1mm,anchor=south] (0,0) node{\small{$-17$}};
	\end{scope}
	\begin{scope}[xshift=15mm]
	\draw[very thick](0,0)--(0,1);
	\draw[xshift=5mm,yshift=-1mm,anchor=south] (0,0) node{\small{$-30$}};
	\end{scope}
	\begin{scope}[xshift=30mm]
	\draw[very thick](0,0)--(1,0);
	\draw[xshift=5mm,yshift=-1mm,anchor=south] (0,0) node{\small{$-18$}};
	\end{scope}
\end{scope}
\begin{scope}[xshift=13.5cm,yshift=17mm]
	\draw[very thick](0,0)--(1,0);
	\draw[xshift=5mm,yshift=-1mm,anchor=south] (0,0) node{\small{$-6$}};
\end{scope}
\begin{scope}[xshift=17cm,yshift=-1cm]
	\draw[gray]
		(0,1)--(1,1)--(1,8)--(0,8)--(0,0)
		(0,2)--(2,2)--(2,8)--(0,8)
		(0,3)--(3,3)--(3,8)--(0,8)
		(0,4)--(4,4)--(4,8)--(0,8)
		(0,5)--(5,5)--(5,8)--(0,8)
		(0,6)--(6,6)--(6,7)--(0,7);
	\draw[very thick] (0,0)--(0,2)--(1,2)--(1,3)--(2,3)--(2,5)--(3,5)--(3,6)--(4,6)--(4,7)--(5,7);
	\draw[xshift=-1mm,yshift=5mm,anchor=west]
		(0,0) node{\small{$-35$}}
		(0,1) node{\small{$-30$}}
		(1,2) node{\small{$-22$}}
		(2,3) node{\small{$-17$}}
		(2,4) node{\small{$-13$}}
		(3,5) node{\small{$-9$}}
		(4,6) node{\small{$-4$}};
	\draw[xshift=5mm,yshift=-1mm,anchor=south]
		(0,2) node{\small{$-23$}}
		(1,3) node{\small{$-18$}}
		(2,5) node{\small{$-11$}}
		(3,6) node{\small{$-6$}}
		(4,7) node{\small{$-1$}};
	\draw[->,thick] (-1.6,4)--(-.8,4);
\end{scope}
\end{tikzpicture}
\caption{The labelling of the steps of a path $\pi$ (left), the set $X$ of labelled steps (middle), and the path $\op{sw}(\pi)$ of steps in increasing order (right).}
\label{Figure:sweepC}
\end{center}
\end{figure}

\subsection{The sweep map}

A generalisation of the zeta map to rational Dyck paths called the sweep map was defined by
\cite{ALW2014SweepMaps}. The concept of the sweep map is as follows. Given a path one assigns to each step a label, the labels being distinct integers. To obtain the image of a path under the sweep map, one rearranges the steps such that the labels are in increasing order.

We now give a description of the zeta map of type $C$ similar to the sweep map on Dyck paths. Given a path $\pi=s_1s_2,\dots,s_{2n}\in\mathcal L_{n,n}$ assign a label to each step by setting $\ell(s_1)=0$, $\ell(s_{i+1})=\ell(s_i)+2n+1$ if $s_i=N$, and $\ell(s_{i+1})=\ell(s_i)-2n$ if $s_i=E$. Now define a collection $X$ of labelled steps as follows. If $\ell(s_i)<0$ then add $(s_i,\ell(s_i))$. If $\ell(s_i)>0$ then add $(s_{i-1},-\ell(s_i))$. Finally, for the step $s_1$ which is the only step labelled $0$, add $(s_{2n},-n)$. Thus, $X$ contains $2n$ labelled steps.

Finally, draw a path as follows. Choose $(s,\ell)\in X$ such that $\ell$ is the minimal label among all pairs in $X$. Draw the step $s$ and remove $(s,\ell)$ from $X$. Repeat until $X$ is empty. We denote the path obtained in this way by $\op{sw}(\pi)$. See \reff{sweepC}. We conclude with the following theorem.

\begin{mythrm}\label{thrm:sweepC} For each lattice path $\pi\in\mathcal L_{n,n}$ we have $\op{sw}(\pi)=\zeta(\pi)$. In particular, the sweep map $\op{sw}:\mathcal L_{n,n}\to\mathcal B_n$ is a bijection.
\end{mythrm}



\nocite{}
\bibliographystyle{abbrvnat}
\bibliography{submission}
\label{sec:biblio}

\end{document}